\documentclass{article}
\usepackage{amsmath,amsthm,amscd,amssymb,amsfonts}
\input{amssym.def}
\input{amssym.tex}

\newtheorem{thm}{Theorem}[section]
\newtheorem{prop}{Proposition}[section]
\newtheorem{lm}{Lemma}[section]

\newtheorem{cor}{Corollary}[section]

\title{An analog of the Iwasawa conjecture for a compact hyperbolic threefold}
\author{Ken-ichi SUGIYAMA
}
\begin{document}
\maketitle
\begin{abstract}
For a local system on a compact hyperbolic threefold, under a
 cohomological assumption, we will show that the
 order of its twisted Alexander polynomial and of the Ruelle L
 function at $s=0$ coincide. Moreover we will show that their leading
 constant are also identical. These results may be considered as a
 solution of a geomeric analogue of the Iwasawa conjecture in the algebraic number theory.
\footnote{1991 Mathematics Subject Classification : 57M25, 57M10, 57N10}
\footnote{2000 Mathematics Subject Classification : 11F32, 11M36, 57M25, 57M27} 
\end{abstract}
\section{Introduction}
In rescent days, it has been recognized there are many similarities
between the theory of a number field and one of a topological
threefold. In this note, we will show one more evidence, which is ``a
geometric analog of the Iwasawa conjecture''.\\

At first let us recall the original Iwasawa conjecture (\cite{Washington}). Let $p$ be an
odd prime and $K_n$ a cyclotomic field ${\mathbb Q}(\zeta_{p^n})$. The
Galois group ${\rm Gal}(K_n/{\mathbb Q})$ which is
isomorphic to ${\mathbb Z}/(p^{n-1})\times {\mathbb F}^{*}_p$ by the
cyclotomic character $\omega$ acts on
the $p$-primary part of the ideal class group $A_n$ of $K_n$. By the
action of ${\rm Gal}(K_1/{\mathbb Q})\simeq {\mathbb F}^{*}_p$, it has a
decomposition
\[
 A_{n}=\oplus_{i=0}^{p-2}A_{n}^{\omega^{i}},
\]
where we set
\[
 A_n^{\omega^i}=\{\alpha\in A_n \,|\, \gamma
      \alpha=\omega(\gamma)^i\alpha \,\,\mbox{for}\,\, \gamma \in {\rm
      Gal}(K_1/{\mathbb Q})\}.
\]
For each $i$ let us take the inverse limit with respect to the norm map:
\[
 X_i=\lim_{\leftarrow}A_n^{\omega^i}.
\]
If we set $K_{\infty}=\cup_{n}K_n$ and $\Gamma={\rm
Gal}(K_{\infty}/K_1)$, each $X_{i}$ becomes a ${\mathbb
Z}_{p}[[\Gamma]]$-module. Since there is an (noncanocal) isomorphism
${\mathbb Z}_{p}[[\Gamma]]\simeq {\mathbb Z}_{p}[[s]]$, each $X_i$ may be considered
as a ${\mathbb Z}_{p}[[s]]$-module. Iwasawa has shown that it is a
torsion ${\mathbb Z}_{p}[[s]]$-module and let ${\mathcal L}_{p}^{alg,i}$
be its generator, which will be referred as {\it the Iwasawa power
series.}\\

On the other hand, let 
\[
 {\mathbb Z}_{p}[[s]]\simeq {\mathbb Z}_{p}[[\Gamma]]\stackrel{\chi}\to{\mathbb Z}_{p}
\]
be the ring homomorphism induced by $\omega.$ For each $0 < i <p-1$, using the Kummer
congruence of the Bernoulli numbers, Kubota-Leopoldt and Iwasawa have
independently constructed an element of  ${\mathcal L}_p^{ana,i}$ which
satisfies
\[
 \chi^{r}({\mathcal L}_p^{ana,i})=(1-p^{r})\zeta(-r),
\]
for any positive integer $r$ which is congruent $i$ modulo $p-1$. Here
$\zeta$ is {\it the Riemann zeta function.} We will refer
${\mathcal L}_p^{ana,i}$ as {\it the $p$-adic zeta function}. The
Iwasawa main conjecture, which has been solved by Mazur and Wiles (\cite{MW})
says that ideals in ${\mathbb Z}_{p}[[s]]$ generated by ${\mathcal L}_{p}^{alg,i}$ and ${\mathcal
L}_p^{ana,i}$ are equal.\\

Now we will explain our geometric analog of the Iwasawa main
conjecture. \\

It is broadly recognized a geometric substitute for the
Iwasawa power series is the Alexander invariant. Let $X$ be a connected
finite CW-complex of dimension three and $\Gamma_{g}$ its fundamental
group. In what follows, we always assume that there is a surjective
homomorphism
\[
 \Gamma_{g} \stackrel{\epsilon}\to {\mathbb Z}.
\]
Let $X_{\infty}$ be the infinite cyclic covering of $X$ which
corresponds to ${\rm Ker}\,\epsilon$ by the geometric Galois theory and
$\rho$ a finite dimensional unitary representation of $\Gamma_g$. Then
$H_{\cdot}(X_{\infty},\,{\mathbb C})$ and $H_{\cdot}(X_{\infty},\,\rho)$
have an action of ${\rm
Gal}(X_\infty/X)\simeq {\mathbb Z}$, which make them $\Lambda$-modules.
Here we set $\Lambda={\mathbb C}[{\mathbb Z}]$ which is isomorphic to the Laurent polynomial ring ${\mathbb
C}[t,\,t^{-1}]$. Suppose that each of them is a torsion
$\Lambda$-module. Then due to the results of Milnor (\cite{MilnorI}), we
know $H^{i}(X_{\infty},\,\rho)$ is also a torsion $\Lambda$-module for
all $i$ and
vanishes for $i\geq 3$. Let $\tau^{*}$ be the action of $t$ on
$H^{i}(X_{\infty},\,\rho)$. Then {\it the Alexander invariant} is
defined to be the alternating product of the characteristic polynomials:
\[
 A_{\rho}^{*}(t)=\frac{\det[t-\tau^{*}\,|\, H^{0}(X_{\infty},\,\rho)]\cdot \det[t-\tau^{*}\,|\, H^{2}(X_{\infty},\,\rho)]}{\det[t-\tau^{*}\,|\, H^{1}(X_{\infty},\,\rho)]}.
\]

On the other hand, we will take {\it the Ruelle L-function} as a geometric substitute for
the $p$-adic zeta function. Let $X$ be a connected closed hyperbolic
threefold. Then its fundamental group $\Gamma_g$ may be considered as a
torsion-free cocompact discrete subgroup of $PSL_2({\mathbb C})$. By the one to one
correspondence between the set of loxiodromic conjugacy classes of
$\Gamma_g$ and one of closed geodesics of $X$,  {\it the Ruelle
L-function} is defined to be a product of the characteristic polynomials
of $\rho(\gamma)$ over prime closed geodesics:
\[
 R_{\rho}(s)=\prod_{\gamma}\det[1-\rho(\gamma)e^{-sl(\gamma)}].
\]
Here $s$ is a complex number and $l(\gamma)$ is the length of
$\gamma$. It absolutely convergents for $s$ whose real part is
sufficiently large. Fried (\cite{Fried}) has shown that it is
meromorphically continued in the whole plane. Using his results
(\cite{Fried}), we will show the following theorem.
\begin{thm}
Suppose that $H^{0}(X_{\infty},\,\rho)$ vanishes and let $\beta$ be the
 dimension of $H^{1}(X,\,\rho)$. Then
\[
 -2\beta= {\rm ord}_{s=0}R_{\rho}(s)\geq 2{\rm ord}_{t=1}{A_{\rho}^{*}(s)},
\]
and the identity holds if the action of $\tau^{*}$ on
 $H^{1}(X_{\infty},\,\rho)$ is semisimple.
Moreover if all $H^{i}(X,\,\rho)$ vanish, we have
\[
 |R_{\rho}(0)|=\delta_{\rho}|A_{\rho}^{*}(1)|^{2},
\]
where $\delta_{\rho}$ is a positive constant which can be determined explicitly. 
\end{thm}
In particular if we make a change of variables:
\[
 t=s+1,
\]
our theorem implies that two ideals in ${\mathbb C}[[s]]$ which are generated by
$R_{\rho}(s)$ and $A_{\rho}^{*}(s)^{2}$ coincide. Thus it may be
considered as a solution of {\it a geometric analog of the Iwasawa main
conjecture.}\\

In fact, after a certain modification, {\bf Theorem 1.1} is still true for a complete hyperbolic threefold of a
finite volume, which will be discussed in \cite{Sugiyama}.\\

When $X_{\infty}$ is homeomorhic to a mapping torus derived a
homeomorphism $f$
of a compace Riemannian surface $S$, we can prove {\it a limit formula}.
\begin{thm}
Suppose that $X$ is homeomorphic to a mapping torus of $(S,\, f)$
 and that the surjective homomorphism $\epsilon$ is induced from the structure map:
\[
 X\to S^{1}.
\]
If $H^{0}(S,\,\rho)$ vanishes, we have
\[
 -2\beta= {\rm ord}_{s=0}R_{\rho}(s)\geq 2{\rm ord}_{t=1}{A_{\rho}^{*}(s)},
\]
and the identiy holds if the action of $f^{*}$ on
 $H^{1}(S,\,\rho)$ is semisimple. Moreover if this is satisfied, we have
\[
 \lim_{s\to 0}|s^{2\beta}R_{\rho}(s)|=\lim_{t\to 1}|(t-1)^{\beta}A_{\rho}^{*}(t)|^{2}=|\tau^{*}_{\mathbb C}(X,\,\rho)|^{2},
\]
where $\tau^{*}_{\mathbb C}(X,\,\rho)$ is the (cohomological)
 Milnor-Reidemeister torsion of $X$ and $\rho$.
\end{thm}

\section{The Milnor-Reidemeister torsion and the Alexander invariant}
Let $\Lambda={\mathbb C}[t,\,t^{-1}]$ be a Laurent polynomial ring of
complex coefficients. The following lemma is easy to see.
\begin{lm}
Let $f$ and $g$ be elements of $\Lambda$ such that
\[
 f=ug,
\]
where $u$ is a unit.
Then their order at $t=1$ are equal:
\[
 {\rm ord}_{t=1}f={\rm ord}_{t=1}g.
\]
\end{lm}
We will recall {\it the  Milnor-Reidemeister torsion} of a complex
(\cite{MilnorW},\cite{MilnorI}).\\

 Let $(C_{\cdot},\,\partial_{\cdot})$ be
a bounded complex of free $\Lambda$-modules of finite rank whose homology
groups are torsion $\Lambda$-modules. Suppose that it is given a base ${\bf c}_i$
for each $C_{i}$.  
Such a complex will
refered as {\it a based complex.} We set
\[
 C_{even}=\oplus_{i\equiv 0(2)} C_{i},\quad C_{odd}=\oplus_{i\equiv 1(2)} C_{i},
\]
which are free $\Lambda$-modules of finite rank with basis ${\bf
c}_{even}=\oplus_{i\equiv 0(2)}{\bf c}_i$ and ${\bf c}_{odd}=\oplus_{i\equiv 1(2)}{\bf c}_i$ respectively. Choose a
base ${\bf b}_{even}$ of a $\Lambda$-submodule $B_{even}$ of $C_{even}$
(necessary free) which is the image of the
differential and column vectors ${\bf x}_{odd}$ of $C_{odd}$
so that
\[
 \partial {\bf x}_{odd}={\bf b}_{even}.
\]  
Similarly we take ${\bf b}_{odd}$ and ${\bf x}_{even}$ satisfying
\[
 \partial {\bf x}_{even}={\bf b}_{odd}.
\] 
Then ${\bf x}_{even}$ and ${\bf b}_{even}$ are expressed by a linear
combination of ${\bf c}_{even}$:
\[
 {\bf x}_{even}=X_{even}{\bf c}_{even},\quad {\bf b}_{even}=Y_{even}{\bf c}_{even},
\]
and we obtain a square matrix
\[
 \left(
\begin{array}{c}
X_{even}\\
Y_{even}
\end{array}
\right).
\] 
Similarly the equation
\[
 {\bf x}_{odd}=X_{odd}{\bf c}_{odd},\quad {\bf b}_{odd}=Y_{odd}{\bf c}_{odd}
\]
yields a square matrix
\[
 \left(
\begin{array}{c}
X_{odd}\\
Y_{odd}
\end{array}
\right). 
\]
Now {\it the Milnor-Reidemeister torsion} $\tau_{\Lambda}(C_{\cdot},{\bf
c}_{\cdot})$ of the based complex $\{C_{\cdot},{\bf c}_{\cdot}\}$ is defined as
\begin{equation}
\tau_{\Lambda}(C_{\cdot},{\bf c}_{\cdot})=
\pm
\frac{\det  \left(
\begin{array}{c}
X_{even}\\
Y_{even}
\end{array}
\right)}
{\det  \left(
\begin{array}{c}
X_{odd}\\
Y_{odd}
\end{array}
\right)}
\end{equation}
It is known $\tau_{\Lambda}(C_{\cdot},{\bf c}_{\cdot})$ is independent
of a choice of ${\bf b}_{\cdot}$.\\

Since $H_{\cdot}(C_{\cdot})$ are torsion $\Lambda$-modules, they are
finite dimensional complex vector spaces. Let $\tau_{i*}$ be the action
of $t$ on $H_{i}(C_{\cdot})$. Then {\it the Alexander invariant}
is defined to be the alternating product of their characteristic polynomials:
\begin{equation}
 A_{C_{\cdot}}(t)=\prod_{i}\det[t-\tau_{i*}]^{(-1)^{i}}.
\end{equation}
Then {\bf Assertion 7} of \cite{MilnorI} shows the fractional
ideals generated by $\tau_{\Lambda}(C_{\cdot},{\bf c}_{\cdot})$ and
$A_{C_{\cdot}}(t)$ are equal:
\[
 (\tau_{\Lambda}(C_{\cdot},{\bf c}_{\cdot}))=(A_{C_{\cdot}}(t)).
\]
In particular {\bf Lemma 2.1} implies
\begin{equation}
{\rm ord}_{t=1}\tau_{\Lambda}(C_{\cdot},{\bf c}_{\cdot})={\rm ord}_{t=1}A_{C_{\cdot}}(t),
\end{equation}
and we know 
\[
\tau_{\Lambda}(C_{\cdot},{\bf c}_{\cdot})=\delta\cdot t^{k}A_{C_{\cdot}}(t), 
\]
where $\delta$ is a non-zero complex number and $k$ is an
integer. $\delta$ will be referred as {\it the difference} of the
Alexander invariant and the Milnor-Reidemeister torsion.\\

Let $\{\overline{C_{\cdot}},\,\overline{\partial}\}$ be a bounded
complex of a finite dimensional vector spaces over ${\mathbb C}$. 
If it is given basis ${\bf c}_i$ and ${\bf h}_i$ for each $\overline{C_{i}}$
 and $H_{i}(\overline{C_{\cdot}})$ respectively, the Milnor-Reidemeister
 torsion 
$\tau_{\mathbb C}(\overline{C_{\cdot}},\,\overline{{\bf c}_{\cdot}})$ 
is also defined (\cite{MilnorW}). Such a complex will be referred as {\it a
 based complex} again. By definition, if the complex is acyclic, it 
 coincides with (1). Let $(C_{\cdot},\,{\bf c}_{\cdot})$ be a based
 bounded complex over $\Lambda$ whose homology groups are
 torsion $\Lambda$-modules. Suppose its annihilator $Ann_{\Lambda}(H_{i}(C_{\cdot}))$ does
 not contain $t-1$ for each $i$. Then 
\[
 (\overline{C_{\cdot}},\,\overline{\partial})=(C_{\cdot},\,{\bf c}_{\cdot})\otimes_{\Lambda}\Lambda/(t-1)
\]
is a based acyclic complex over ${\mathbb C}$ with a preferred base $\overline{{\bf
c}_{\cdot}}$ which is the reduction of ${\bf c}_{\cdot}$ modulo $(t-1)$. This observation shows the
following proposition.
\begin{prop}
Let $(C_{\cdot},\,{\bf c}_{\cdot})$ be a based
 bounded complex over $\Lambda$ whose homology groups are
 torsion $\Lambda$-modules. Suppose the annihilator $Ann_{\Lambda}(H_{i}(C_{\cdot}))$ does
 not contain $t-1$ for each $i$. Then we have
\[
 \tau_{\Lambda}(C_{\cdot},\,{\bf c}_{\cdot})|_{t=1}=\tau_{\mathbb C}(\overline{C_{\cdot}},\,\overline{{\bf c}_{\cdot}})
\]

\end{prop}

For a later purpose we will consider these dual.\\

 Let
$\{C^{\cdot},\,d\}$ be the dual complex of $\{C_{\cdot},\,\partial\}$:
\[
 (C^{\cdot},\,d)=Hom_{\Lambda}((C_{\cdot},\,\partial),\,\Lambda).
\]
By the universal coefficient theorem we have
\[
 H^{q}(C^{\cdot},\,d)=Ext^{1}_{\Lambda}(H_{q-1}(C_{\cdot},\,\partial),\,\Lambda)
\]
and the cohomology groups are torsion $\Lambda$-modules. Moreover the
characteristic polynomial of $H^{q}(C^{\cdot},\,d)$ is equal to
one of $H_{q-1}(C_{\cdot},\,\partial)$. Thus if we define {\it the
Alexander invariant} $A_{C^{\cdot}}(t)$ of $\{C^{\cdot},\,d\}$ by the
same way as (2),
we have
\begin{equation}
 A_{C^{\cdot}}(t)=A_{C_{\cdot}}(t)^{-1}.
\end{equation}

\section{The Milnor-Reidemeister torsion of a CW-complex of dimension
 three}
Let $X$ be a connected finite CW-complex and $\{c_{i,\alpha}\}_\alpha$
its $i$-dimensional cells. We will fix its base point $x_0$ and let
$\Gamma$ be the fundamental group of $X$ with a base point $x_{0}$. Let
$\rho$ be a unitary representation of a finite dimension and $V_{\rho}$
its representation space. Suppose that there is a surjective
homomorphism
\[
 \Gamma \stackrel{\epsilon}\to {\mathbb Z},
\]
and let $X_{\infty}$ be the infinite cyclic covering of $X$ which
corresponds to ${\rm Ker}\,\epsilon$ by the Galois theory. Finally let
$\tilde{X}$ be the universal covering of $X$. \\

The chain complex $(C_{\cdot}(\tilde{X}),\,\partial)$ is a complex of
free ${\bf C}[\Gamma]$-module of finite rank. We take a
lift of ${\bf c}_i=\{c_{i,\alpha}\}_\alpha$ as a base of
$C_{i}(\tilde{X})$, which will be also denoted by the same character. Note that such a choice of base has an ambiguity of
the action of $\Gamma$.\\

Following \cite{KL} consider a complex over ${\mathbb C}$:
\[
 C_{i}(X,\,\rho)=C_{i}(\tilde{X})\otimes_{{\mathbb C}[\Gamma]} V_{\rho}.
\]

On the other hand, restricting $\rho$ to ${\rm Ker}\,\epsilon$, we will
make a chain complex
\[
 C_{\cdot}(X_{\infty},\,\rho)=C_{\cdot}(\tilde{X})\otimes_{{\mathbb C}[{\rm Ker}\epsilon]}V_{\rho},
\]
which has the
following description. Let us consider ${\mathbb C}[{\mathbb
Z}]\otimes_{\mathbb C}V_{\rho}$ as $\Gamma$-module by
\[
 \gamma(p\otimes v)=p\cdot t^{\epsilon(\gamma)}\otimes \rho(\gamma)\cdot
 v,\quad p\in {\mathbb C}[{\mathbb Z}],\,v\in V_{\rho}.
\]
Then $C_{\cdot}(X_{\infty},\,\rho)$ is isomorphic to a complex
(\cite{KL} {\bf Theorem 2.1}):
\[
 C_{\cdot}(X,\,V_{\rho}[{\mathbb
 Z}])=C_{\cdot}(\tilde{X})\otimes_{\mathbb{C}[\Gamma]}({\mathbb
 C}[{\mathbb Z}]\otimes_{\mathbb C}V_{\rho}).
\]
and we know $C_{\cdot}(X_{\infty},\,\rho)$ is a bounded complex of free
$\Lambda$-modules of finite rank. We will fix a unitary base ${\bf v}=\{v_1,\cdots,v_m\}$ of
 $V_{\rho}$ and make it a based complex with a preferred base ${\bf
c}_{\cdot}\otimes{\bf v}=\{c_{i,\alpha}\otimes v_j\}_{\alpha,i,j}$.\\

 In the following we will fix an isomorphism between ${\mathbb
C}[{\mathbb Z}]$ and $\Lambda$ and identify them. Note that such an
isomorphism is determined modulo $t^{k}\,(k\in {\mathbb Z})$. Note that
for such a choice there is an ambiguity of sending the generator $1$ to
$t^{\pm 1}$. Then by the surjection:
\[
 \Lambda \to \Lambda/(t-1)\simeq {\mathbb C},
\]
$C_{\cdot}(X_{\infty},\,\rho)\otimes_{\Lambda}{\mathbb C}$ is isomorphic to
$C_{\cdot}(X,\,\rho)$. Moreover if we take ${\bf
c}_{\cdot}\otimes{\bf v}$ as a base of the latter, they are
isomorphic as based complexes.\\

Let $C^{\cdot}(\tilde{X})$ be the cochain complex of $\tilde{X}$:
\[
 C^{\cdot}(\tilde{X})=Hom_{{\mathbb C}[\Gamma]}(C_{\cdot}(\tilde{X}),\,{\mathbb C}[\Gamma]),
\]
which is a bounded complex of free ${\mathbb C}[\Gamma]$-module of a finite rank. For
each $i$ we will take the dual ${\bf c}^i=\{c^{i}_{\alpha}\}_\alpha$ of
${\bf c}_i=\{c_{i,\alpha}\}_\alpha$ as a base of
$C^{i}(\tilde{X})$. Thus $C^{\cdot}(\tilde{X})$ becomes a based complex
with a preferred base ${\bf c}^{\cdot}=\{{\bf c}^{i}\}_i$. Since $\rho$
is a unitary representation, it is easy to
see that the dual complex of $C_{\cdot}(X_{\infty},\,\rho)$ is
isomorphic to 
\[
 C^{\cdot}(X_{\infty},\,\rho)=C^{\cdot}(\tilde{X})\otimes_{{\mathbb
 C}[\Gamma]}(\Lambda\otimes_{{\mathbb C}}V_{\rho}),
\]
if we twist its complex structure by the complex
conjugation. Also we will make it a based complex by the base ${\bf
c}^{\cdot}\otimes {\bf v}=\{c^{i}_{\alpha}\otimes v_j\}_{\alpha,i,j}$.\\

Dualizing the exact sequence 
\[
 0 \to C_{\cdot}(X_{\infty},\,\rho) \stackrel{t-1}\to
 C_{\cdot}(X_{\infty},\,\rho) \to C_{\cdot}(X,\,\rho) \to 0
\] 
in the derived category of bounded complex of finitely generated
$\Lambda$-modules, we will obtain a distinguished triangle:
\begin{equation}
  C^{\cdot}(X,\,\rho) \to C^{\cdot}(X_{\infty},\,\rho)\stackrel{t-1}\to
 C^{\cdot}(X_{\infty},\,\rho) \to  C^{\cdot}(X,\,\rho)[1]\to .
\end{equation} 
Here we set 
\[
 C^{\cdot}(X,\,\rho)=C^{\cdot}(\tilde{X},\,\rho)\otimes_{{\mathbb C}[\Gamma]}V_{\rho}.
\]
and in general for a bounded complex $C^{\cdot}$, $C^{\cdot}[n]$ denotes its
{\it shift}, which is defined as
\[
 C^{i}[n]=C^{i+n}.
\]
Note that $C^{\cdot}(X,\,\rho)$ is isomorphic to the reduction of 
$C^{\cdot}(X_{\infty},\,\rho)$ modulo $(t-1)$.\\

Let $\tau^{*}$ be the action of $t$ on
$H^{\cdot}(X_{\infty},\,\rho)$. Then (5) induces an exact sequence:
\begin{equation}
 \to H^{q}(X,\,\rho)\to H^{q}(X_{\infty},\,\rho)\stackrel{\tau^{*}-1}\to
 H^{q}(X_{\infty},\,\rho) \to H^{q+1}(X,\,\rho)\to. 
\end{equation} 
In the following, we will assume that the dimension of $X$ is three and
that all $H_{\cdot}(X_{\infty},\,{\mathbb C})$ and
$H_{\cdot}(X_{\infty},\,\rho)$ are finite dimensional vector spaces over
${\mathbb C}$. The arguments of $\S 4$ of \cite{MilnorI} will show 
the following theorem.
\begin{thm}(\cite{MilnorI}) 
\begin{enumerate}
\item 
For $i\geq 3$, 
 $H^{i}(X_{\infty},\,\rho)$ vanishes. 
\item For $0\leq i \leq 2$, $H^{i}(X_{\infty},\,\rho)$ is a finite dimensional vector
      space over ${\mathbb C}$ and there is a perfect pairng:
\[
 H^{i}(X_{\infty},\,\rho)\times H^{2-i}(X_{\infty},\,\rho)\to {\mathbb C}.
\] 
\end{enumerate}
\end{thm}
 The perfect pairing will be referred as {\it the Milnor duality}.\\

 Let
 $A_{\rho*}(t)$ and $A^{*}_{\rho}(t)$ be the Alexander invariants of
 $C_{\cdot}(X_{\infty},\,\rho)$ and $C^{\cdot}(X_{\infty},\,\rho)$
 respectively. Since the latter complex is the dual of the previous one,
 (4) implies
\[
 A^{*}_{\rho}(t)=A_{\rho*}(t)^{-1}.
\]
Let $\tau^{*}_{\Lambda}(X_{\infty},\,\rho)$ be the Milnor-Reidemeister
torsion of $C^{\cdot}(X_{\infty},\,\rho)$ with respect to the preferred
base ${\bf c}^{\cdot}\otimes{\bf v}$. Because of an ambiguity of a
choice of ${\bf c}^{\cdot}$ and ${\bf v}$, it is well-defined modulo
\[
 \{zt^{n}\,|\,z\in {\mathbb C},\,|z|=1,\,n\in{\mathbb Z}\}.
\]
Let $\delta_{\rho}$ be the absolute value of the difference between
$A^{*}_{\rho}(t)$ and $\tau^{*}_{\Lambda}(X_{\infty},\,\rho)$. The discussion of the
previous section shows the following theorem.
\begin{thm}
The order of $\tau^{*}_{\Lambda}(X_{\infty},\,\rho)$, $A^{*}_{\rho}(t)$
 and $A_{\rho*}(t)^{-1}$ at $t=1$ are equal. Let $\beta$ be the
 order. Then we have
\begin{eqnarray*}
 \lim_{t\to 1}|(t-1)^{-\beta}\tau^{*}_{\Lambda}(X_{\infty},\,\rho)|&=&\delta_{\rho}\lim_{t\to 1}|(t-1)^{-\beta}A^{*}_{\rho}(t)|\\
&=& \delta_{\rho}\lim_{t\to 1}|(t-1)^{-\beta}A_{\rho*}(t)^{-1}|.
\end{eqnarray*}
\end{thm}
Note that {\bf Theorem 3.1} implies that the Alexander invariant has the
following form:
\begin{equation}
 A^{*}_{\rho}(t)=\frac{\det[t-\tau^{*}\,|\,H^{0}(X_{\infty},\,\rho)]\cdot
 \det[t-\tau^{*}\,|\,H^{2}(X_{\infty},\,\rho)]}{\det[t-\tau^{*}\,|\,H^{1}(X_{\infty},\,\rho)]}.
\end{equation}
\begin{thm}
Suppose $H^{0}(X_{\infty},\,\rho)$ vanishes. Then we have
\[
 {\rm ord}_{t=1} A^{*}_{\rho}(t)\leq -\dim H^{1}(X,\,\rho),
\]
and the identity holds if the action of $\tau^{*}$ on
 $H^{1}(X_{\infty},\,\rho)$ is semisimple.
\end{thm}
{\bf Proof.} We know by {\bf Theorem 3.1} that
$H^{2}(X_{\infty},\,\rho)$ also vanishes. Moreover the exact sequence (6) shows
\[
 0 \to H^{1}(X,\,\rho) \to H^{1}(X_{\infty},\,\rho)
 \stackrel{\tau^{*}-1} \to H^{1}(X_{\infty},\,\rho),
\] 
which implies
\[
 {\rm ord}_{t=1}\det[t-\tau^{*}\,|\,H^{1}(X_{\infty},\,\rho)]\geq\dim H^{1}(X,\,\rho),
\]
and the identity holds if $\tau^{*}$ is semisimple.
Now the desired result follows from (7).
\begin{flushright}
$\Box$
\end{flushright}
Let $\tau^{*}_{{\mathbb C}}(X,\,\rho)$ be the Milnor-Reidemeister
torsion of $C^{\cdot}(X,\,{\rho})$ with respect to the preferred base
${\bf c}^{\cdot}\otimes{\bf v}$. By an ambiguity of a choice of ${\bf v}$, only its
absolute value is well-defined. 
\begin{thm}
Suppose $H^{i}(X,\,\rho)$ vanishes for all $i$. Then we have
\[
 |\tau^{*}_{{\mathbb C}}(X,\,\rho)|=\delta_{\rho}|A^{*}_{\rho}(1)|=\frac{\delta_{\rho}}{|A_{\rho*}(1)|}.
\]
\end{thm}
{\bf Proof.} The exact sequence (6) and the assumption implies $t-1$ is
not contained in the annihilator of $H^{\cdot}(X_{\infty},\,\rho)$. Now
the theorem will follow from {\bf Proposition 2.1} and {\bf Theorem 3.2}.
\begin{flushright}
$\Box$
\end{flushright}
In the following, we will specify these arguments to a mapping torus.\\

 Let $S$ be a
connected CW-complex of dimension 2 and $f$ its automorphism. Let $X$ be
the mapping torus of the pair $S$ and $f$. We will take a base point from
$S$ and let $\Gamma_S$ be the fundamental group of $S$ with respect to
the point. Supoose $H^{0}(S,\,\rho)$ vanishes.
(e.g. The restriction of $\rho$ to $\Gamma_S$ is irreducible.) Let 
\[
 \Gamma = \pi_1(X,\,s_0)\stackrel{\epsilon}\to {\mathbb Z}
\]
be the homomorphism induced by the structure map:
\[
 X \to S^{1}.
\]
Then $X_{\infty}$ is a product of $S$ with the real axis and therefore
$C^{\cdot}(X_{\infty},\,\rho)$ is chain homotopic to
$C^{\cdot}(S,\,\rho)$.  We have an exact sequence of complexes:
\begin{equation}
 0 \to C^{\cdot}(S,\,\rho)[-1] \to C^{\cdot}(S,\,\rho) \to
 C^{\cdot}(X,\,\rho) \to 0.
\end{equation}
The cells of $S$ defines a base ${\bf c}_S$ of the chain complex
$C_{\cdot}(S)$. 
Thier product with the unit interval and themselves form a base
$C_{\cdot}(X)$, which will be denoted by ${\bf c}_X$. Let
${\bf c}^{S}$ and ${\bf c}^{X}$ be the dual of them. If we fix a unitary
base ${\bf v}=\{v_1,\cdots,v_m\}$, both $C^{\cdot}(S,\,\rho)$ and
$C^{\cdot}(X,\,\rho)$ become based complexes with preferred base ${\bf
c}^{S}\otimes{\bf v}$ and ${\bf c}^{X}\otimes{\bf v}$ respectively. \\

The exact sequence of cohomology groups of (8) may be considered as an
acyclic complex:
\begin{equation}
\begin{array}{cccccc}
 &0&\to&H^{0}(X,\,\rho)&\to&H^{0}(S,\,\rho)\\
\stackrel{f^{*}-1}\to&H^{0}(S,\,\rho)&\to&H^{1}(X,\,\rho)&\to&H^{1}(S,\,\rho)\\
\stackrel{f^{*}-1}\to&H^{1}(S,\,\rho)&\to&H^{2}(X,\,\rho)&\to&H^{2}(S,\,\rho)\\
\stackrel{f^{*}-1}\to&H^{2}(S,\,\rho)&\to&H^{3}(X,\,\rho)&\to&0,
\end{array}
\end{equation}
which will be denoted by ${\mathcal H}^{\cdot}$. Note that this is
nothing but $(6)$. Since $H^{0}_{c}(S,\,\rho)$ which is a subgroup of $H^{0}(S,\,\rho)$
vanishes, the Poincar\'{e} duality implies $H^{2}(S,\,\rho)$ also does. Thus ${\mathcal H}^{\cdot}$ is isomorphic to 
\[
 {\mathcal H}^{\cdot}_{0}=[H^{1}(X,\,\rho)\to H^{1}(S,\,\rho)
 \stackrel{f^{*}-1}\to H^{1}(S,\,\rho)\to H^{1}(X,\,\rho)][-4].
\]
Choosing basis of $H^{1}(X,\,\rho)$ and $H^{1}(S,\,\rho)$, we will make
${\mathcal H}^{\cdot}_{0}$ a based acyclic complex.\\

Now we compute the Milnor-Reidemeister torsion of
$C^{\cdot}(X,\,\rho)$. In the following computation, we will assume that
the action of $f^{*}$ on $H^{1}(S,\,\rho)$ is semisimple. \\

Since in general the Milnor-Reidemeister torsion
$\tau^{*}_{\mathbb C}(C^{\cdot}[n])$ of a shift of a based bounded complex
$C^{\cdot}$ is equal to $\tau^{*}_{\mathbb C}(C^{\cdot})^{(-1)^n}$, {\bf Theorem 3.2} of \cite{MilnorW} implies
\[
 \tau^{*}_{\mathbb C}(X,\rho)=\tau^{*}_{\mathbb C}({\mathcal H}^{\cdot}_{0}).
\]
In particular we know $\tau^{*}_{\mathbb C}({\mathcal H}^{\cdot}_{0})$ is independent of a choice of
basis of $H^{1}(X,\,\rho)$ and $H^{1}(S,\,\rho)$. We set 
\[
 I=H^{1}(S,\,\rho)/({\rm Ker}[f^{*}-1]),
\]
and let 
\[
 H^{1}(S,\,\rho)\stackrel{\pi}\to I
\]
be the natural projection.
Then $f^{*}-1$ induces an isomorphism of $I$ and we have a diagram:
\begin{equation}
\begin{array}{ccccccccc}
0&\to&H^{1}(X,\,\rho)&\to&H^{1}(S,\,\rho)&\stackrel{\pi}\to&I&\to&0\\
 & &\downarrow& &\downarrow& & \downarrow& & \\
0&\to&H^{1}(X,\,\rho)&\to&H^{1}(S,\,\rho)&\stackrel{\pi}\to&I&\to&0.
\end{array}
\end{equation}
Here the left vertical arrow is the null homomorphism and the middle one is
$f^{*}-1$. The right vertical arrow is the isomorphism induced by it.
The snake lamma yields an exact sequence of acyclic complexes:
\[
 0 \to {\mathcal F}^{\cdot} \to {\mathcal H}^{\cdot}_{0} \to {\mathcal
 G}^{\cdot} \to 0, 
\]
where we set
\[
 {\mathcal F}^{\cdot}=[H^{1}(X,\,\rho)\stackrel{id}\to
 H^{1}(X,\,\rho)\stackrel{0}\to H^{1}(X,\,\rho)\stackrel{id}\to
 H^{1}(X,\,\rho)]
\]
and
\[
 {\mathcal G}^{\cdot}=[0\to I \stackrel{f^{*}-1}\to I \to 0]=[I \stackrel{f^{*}-1}\to I][-1].
\]
We choose basis ${\bf h}$ and ${\bf i}$ of $H^{1}(X,\,\rho)$ and $I$
respectively. Then ${\bf h}$ and a lift of ${\bf i}$ form a base
of $H^{1}(S,\,\rho)$. \\

Now we compute the Milnor-Reidemeister torsion of ${\mathcal
H}^{\cdot}_{0}$.\\

{\bf Theorem 3.1} of \cite{MilnorW} shows
\[
 \tau_{\mathbb C}^{*}({\mathcal H}^{\cdot}_{0})=\tau_{\mathbb
 C}^{*}({\mathcal F}^{\cdot})\cdot \tau_{\mathbb C}^{*}({\mathcal G}^{\cdot}).
\]
Since $\tau_{\mathbb C}^{*}({\mathcal F}^{\cdot})=1$ and since 
\[
 \tau_{\mathbb C}^{*}({\mathcal G}^{\cdot})=\tau_{\mathbb C}^{*}([I
 \stackrel{f^{*}-1}\to I])^{-1}=(\det [f^{*}-1\,|\,I])^{-1},
\]
we have
\begin{equation}
 \tau^{*}_{\mathbb C}(X,\rho)=\tau_{\mathbb C}^{*}({\mathcal H}^{\cdot}_{0})=(\det [f^{*}-1\,|\,I])^{-1}.
\end{equation}
Note that (7) implies
\[
  A^{*}_{\rho}(t)=\frac{1}{\det[t-\tau^{*}\,|\,H^{1}(S,\,\rho)]}. 
\]
Let $\beta$ be the dimension of $H^{1}(X,\,\rho)$. Then {\bf Theorem
3.3} and (11) show that the order of $A^{*}_{\rho}(t)$ is $-\beta$ and that
\[
 \lim_{t\to 1}|(t-1)^{\beta}A^{*}_{\rho}(t)|=|\tau^{*}_{\mathbb C}(X,\rho)|.
\]  
Thus we have proved the following theorem.
\begin{thm}
Let $f$ be an automorphism of a connected finite CW-complex of dimension
 two $S$ and $X$ its mapping torus. Let $\rho$ be a unitary representation
 of the fundamental group of $X$ which satisfies $H^{0}(S,\,{\rho})=0$. Suppose that the surjective
 homomorphism
\[
 \Gamma \stackrel{\epsilon}\to {\mathbb Z}
\]
is induced from the structure map
\[
 X\to S^{1},
\]
and that the action of $f^{*}$ on $H^{1}(S,\,\rho)$ is semisimple.
Then the order of $A^{*}_{\rho}(t)$ is $-\beta$, where $\beta$ be the
 dimension of $H^{1}(X,\,\rho)$ and 
\[
 \lim_{t\to 1}|(t-1)^{\beta}A^{*}_{\rho}(t)|=|\tau^{*}_{\mathbb C}(X,\rho)|.
\]  
\end{thm}
In particular we know that $|\tau^{*}_{\mathbb C}(X,\rho)|$ is
determined by the homotopy class of $f$. Note that without semisimplicity of $f^{*}$, we only have
\[
 {\rm ord}_{t=1}A^{*}_{\rho}(t) \leq -\beta.
\]
\section{Comparison of the Milnor-Reidemeister torsion with the twisted
 Alexander polynomial}
Let $K$ be a knot in the three dimensional sphere and $X$ its
complement. Since the first homology group of $X$ is the infinite
cyclic group, there is a surjective homomorphism $\epsilon$ from the
fundamental group $\Gamma$ of $X$ to ${\mathbb Z}$. Let $\rho$ be an
$m$-dimensional 
unitary representation of $\Gamma$.  We assume the homology groups of
$C_{\cdot}(X_{\infty},\,\rho)$ are
torsion $\Lambda$-modules. Thus $C_{\cdot}(X_{\infty},\,\rho)\otimes_{\Lambda}{\mathbb
C}(t)$ is an acyclic complex over ${\mathbb C}(t)$ and we can define the
twisted Alexander polynomial after Kitano (\cite{Kitano}).\\

Let
\[
 P(\Gamma)=<x_{1},\cdots,x_{k}\,|\,r_{1},\cdots,r_{k-1}>
\]
be the Wirtinger representation of $\Gamma$ and $\gamma$ the  natural
projection from the free group $F_{k}$ with $k$ generators:
\[
 F_{k}\stackrel{\gamma}\to \Gamma.
\]
 The $\rho$ and $\epsilon$
induces a ring homomorphism:
\[
 {\mathbb C}[\Gamma]\stackrel{\rho\otimes \epsilon}\to {\rm
 M}_m({\mathbb C}(t)).
\]
Composing this with $\gamma$, we obtain a homomorphism:
\[
 {\mathbb C}[F_k]\stackrel{\Phi}\to {\rm
 M}_m({\mathbb C}(t)).
\]
Using a CW-complex structure of $X$ which is derived from the
Wirtinger representation, the complex $C_{\cdot}(X_{\infty},\,\rho)\otimes_{\Lambda}{\mathbb
C}(t)$ becomes (\cite{Kitano} p.438)
\begin{equation}
 0\to ({\mathbb C}(t)^{\oplus m})^{\oplus (k-1)}
 \stackrel{\partial_2}\to ({\mathbb C}(t)^{\oplus m})^{\oplus
  k}\stackrel{\partial_1}\to {\mathbb C}(t)^{\oplus m} \to 0.
\end{equation}
Here differentials are 
\[
 \partial_2=
\left(
\begin{array}{ccc}
\Phi(\frac{\partial r_1}{\partial x_1})&\cdots& \Phi(\frac{\partial r_{k-1}}{\partial x_1})\\
\vdots&\ddots&\vdots\\
\Phi(\frac{\partial r_1}{\partial x_{k}})&\cdots& \Phi(\frac{\partial r_{k-1}}{\partial x_k})
\end{array}
\right),
\]
and
\[
 \partial_1=
\left(
\begin{array}{ccc}
\Phi(x_1-1),\cdots,\Phi(x_k-1)
\end{array}
\right),
\]
and the derivatives are taken according to the Fox's free differential
calculus. Note that each entry is an element of ${\rm M}_m({\mathbb
C}(t))$. If we set
\[
 a_{1}=
\left(
\begin{array}{ccc}
\Phi(\frac{\partial r_1}{\partial x_1})&\cdots& \Phi(\frac{\partial r_{k-1}}{\partial x_1})
\end{array}
\right)
\]
and
\[
 A_{1}=
\left(
\begin{array}{ccc}
\Phi(\frac{\partial r_1}{\partial x_2})&\cdots& \Phi(\frac{\partial r_{k-1}}{\partial x_2})\\
\vdots&\ddots&\vdots\\
\Phi(\frac{\partial r_1}{\partial x_{k}})&\cdots& \Phi(\frac{\partial r_{k-1}}{\partial x_k})
\end{array}
\right),
\]
we have
\[
 \partial_2=
\left(
\begin{array}{c}
a_{1}\\
A_{1}
\end{array}
\right).
\]
Kitano(\cite{Kitano} {\bf Proposition 3.1}) has shown that $\Phi(x_j-1)$ and $A_1$ are
invertible and that the
torsion of (12), which is equal to {\it the twisted Alexander
polynomial} $\Delta_{K,\rho}$ of $K$, is
\[
 \frac{\det A_{1}}{\det \Phi(x_1-1)}.
\]
Now let us consider the dual of (12). Since the transpose of
$\partial_2$ is
\[
 \partial_2^{t}=(a_{1}^{t}, A_{1}^{t}),
\]
and since $A_1^{t}$ is invertible, we may take a lift of the standard base
${\bf e}^{*}$ of the dual space of $({\mathbb C}(t)^{\oplus m})^{\oplus
(k-1)}$ as
\[
 \tilde{{\bf e}^{*}}=
\left(
\begin{array}{c}
0\\
(A_1^{t})^{-1}({\bf e}^{*})
\end{array}
\right).
\]
Let ${\bf f}^{*}$ be the standard base of the dual of ${\mathbb
C}(t)^{\oplus m}$. Then the Milnor-Reidemeister torsion of the dual of (12) is
\[
 \det (\partial_1^{t}({\bf f}^{*}),\tilde{{\bf e}^{*}})
=\det
\left(
\begin{array}{cc}
\Phi(x_1-1)^{t}&0\\
\vdots& \\
\Phi(x_k-1)^{t}&(A_1^{t})^{-1}
\end{array}
\right)
=\frac{\det \Phi(x_1-1)}{\det A_1},
\]
which is the inverse of the twisted Alexander polynomial. Thus we have
proved the following theorem.
\begin{thm}
Let $X$ be the complement of a knot $K$ in the three dimensional sphere
 and $\rho$ a unitary representaion of its fundamental group. Suppose
 $H_{i}(X_{\infty},\,\rho)$ are finite dimensional complex vector spaces for all
 $i$. Then the Milnor-Reidemeister torsion of
 $C^{\cdot}(X_{\infty},\,\rho)\otimes_{\Lambda}{\mathbb C}(t)$ is the
 inverse of the twisted Alexander polynomial $\Delta_{K,\rho}$ of $K$. 
\end{thm}
{\bf Proposition 2.1}, {\bf Theorem 3.2}, {\bf Theorem 3.3} and {\bf
Theorem 4.1} imply
\begin{cor}
Suppose $H^{0}(X_{\infty},\,\rho)$ vanishes. Then we have
\[
 {\rm ord}_{t=1}\Delta_{K,\rho}(t)=-{\rm ord}_{t=1}A_{\rho}^{*}(t)\geq \dim H^{1}(X,\,\rho),
\]
and the identity holds if the action of $\tau^{*}$ on
 $H^{1}(X_{\infty},\,\rho)$ is semisimple.
Moreover suppose $H^{i}(X,\,\rho)$ vanishes for all $i$. Then 
\[
 |\tau_{\mathbb C}^{*}(X,\,\rho)|=\frac{1}{|\Delta_{K,\rho}(1)|}.
\]
\end{cor}
\section{A special value of the Ruelle L-function}
Let $X$ be a compact hyperbolic threefold. Then its fundamental group is
identified with a discrete cocompact subgroup $\Gamma$ of
$PSL_2({\mathbb C})$. In particular there is a one to one correspondence
between the set of loxidromic conjugacy classes of $\Gamma$ and the set
of closed geodesics of $X$. Let 
\[
 \Gamma \stackrel{\rho}\to U(V_{\rho})
\]
be a finite dimensional unitary representation. Using the correpondence, for a closed
geodesic $\gamma$ of $X$, we can define a function:
\[
 \det[1-\rho(\gamma)e^{-sl(\gamma)}].
\]
Here $s$ is a complex number and $l(\gamma)$ is the length of $\gamma$. Following \cite{Fried},
we will define {\it the Ruelle L-function} to be
\[
 R_{\rho}(s)=\prod_{\gamma}\det[1-\rho(\gamma)e^{-sl(\gamma)}],
\]
where $\gamma$ runs through prime closed geodesics. 
It is known that $R_{\rho}(s)$ absolutely convergents for ${\rm Re}\,s
>>0$.\\

Let us take a triangulation of $X$ and a unitary base ${\bf
v}=\{v_1,\cdots,v_m\}$ of $V_{\rho}$. Using them, we make a base of the
chain complex $C_{\cdot}(X,\,\rho)$ and its 
Milnor-Reidemeister torsion $\tau^{*}_{\mathbb C}(X,\rho)$ is defined. It is known
that $\tau^{*}_{\mathbb C}(X,\rho)$ is independent of a choice of a
triangulation. Moreover its absolute value is also independent of a
choice a unitary base.\\

If we apply {\bf Theorem 3} of
\cite{Fried} to our case, we will obtain the following theorem.
\begin{thm}
The Ruelle L-function is meromorphically continued in the whole
 plane. Its order at $s=0$ is
\[
 e=4\dim H^{0}(X,\,\rho)-2\dim H^{1}(X,\,\rho).
\]
Moreover we have
\[
 \lim_{s\to 0}|s^{-e}R_{\rho}(s)|=|\tau^{*}_{\mathbb C}(X,\rho)|^2.
\]
\end{thm}
Fried has shown the theorem for an orthogonal representation but his
proof is still valid for a unitary case.\\

Suppose there is a surjective homomorphism 
\[
 \Gamma \stackrel{\epsilon}\to {\mathbb Z},
\]
and let $X_{\infty}$ be the infinite cyclic covering of $X$ which
corresponds to ${\rm Ker}\,{\epsilon}$. Suppose that both
$H_{i}(X_{\infty},\,{\mathbb C})$ and $H_{i}(X_{\infty},\,\rho)$ are finite dimensional
vector spaces for all $i$. 
Then {\bf Theorem 3.3}, {\bf Theorem 3.4} and {\bf Theorem 5.1} implies
\begin{thm}
\begin{enumerate} 
\item Suppose $H^{0}(X,\,\rho)$ vanishes. Then we have
\[
 -2\beta={\rm ord}_{s=0}R_{\rho}(s)\leq 2{\rm ord}_{t=1}A_{\rho}^{*}(t),
\]
where $\beta$ is $\dim H^{1}(X,\,\rho)$. Moreover if the action of
      $\tau^{*}$ on $H_{1}(X_{\infty},\,\rho)$ is semisimple, the
      identity holds.
\item Suppose $H^{i}(X,\,\rho)$ vanishes for all $i$. Then
\[
 |R_{\rho}(0)|=|\delta_{\rho}\cdot A^{*}_{\rho}(1)|^{2}=
\left|
\frac{\delta_{\rho}}{A_{\rho*}(1)}
\right|^2,
\]
where $\delta_{\rho}$ is the difference of the Alexander invariant and the
      Milnor-Reidemeister torsion.
\end{enumerate}
\end{thm}
When $X$ is a mapping torus of an automorphism $f$ of a compact Riemannian
surface $S$, we can say much more. Suppose $\epsilon$ is induced from the
structure map:
\[
 X \to S^{1}.
\]  
Then {\bf Theorem 3.5} implies the following theorem.
\begin{thm} Suppose $H^{0}(S,\,\rho)$ vanishes. Then 
\[
 -2\beta={\rm ord}_{s=0}R_{\rho}(s)\leq2{\rm ord}_{t=1}A_{\rho}^{*}(t),
\]
where $\beta$ is $\dim H^{1}(X,\,\rho)$. Moreover suppose that the action of
 $f^{*}$ on $H^{1}(S,\,\rho)$ is semisimple. Then the identity holds and
 we have 
\[
 \lim_{s\to 0}|s^{2\beta}R_{\rho}(s)|=\lim_{t\to 1}|(t-1)^{\beta}A^{*}_{\rho}(t)|^2=|\tau^{*}_{\mathbb C}(X,\rho)|^2.
\]
\end{thm}
If we make a change of variables:
\[
 t=s+1,
\] 
{\bf Theorem 5.2} and {\bf Theorem 5.3} show the fractional ideals in
the formal power series ring ${\mathbb C}[[s]]$ 
generated by $R_{\rho}(s)$ and $A_{\rho}^{*}(s)^2$ coincide. Thus our
theorems may be considered as a solution of ``a geometric Iwasawa conjecture''.

\vspace{10mm}
\begin{flushright}
Address : Department of Mathematics and Informatics\\
Faculty of Science\\
Chiba University\\
1-33 Yayoi-cho Inage-ku\\
Chiba 263-8522, Japan \\
e-mail address : sugiyama@math.s.chiba-u.ac.jp
\end{flushright}
\end{document}